\documentclass[10pt]{amsart}
\usepackage{amsmath}
\usepackage{amsthm}
\usepackage{amscd}
\usepackage{amssymb}
\usepackage{amsfonts}
\usepackage{graphics}
\usepackage{color}
\usepackage{comment}
\usepackage[dvipsnames]{xcolor}
\usepackage[all,cmtip]{xy}
\usepackage[normalem]{ulem}
\date{}

\newlength{\defbaselineskip }
 \long\def\salta#1{\relax}

 \theoremstyle{plain}
\newtheorem{theorem}{Theorem}[section]
\newtheorem{proposition}[theorem]{Proposition}
\newtheorem{lemma}[theorem]{Lemma}

\theoremstyle{definition}
\newtheorem{definition}[theorem]{Definition}

\newtheorem{remark}[theorem]{Remark}

\newcommand{\cB}{{\mathcal H}}

\newcommand{\weak}{\rightharpoonup}

\newcommand{\na}{\mathbb{N}}

\newcommand{\re}{\mathbb{R}}

\newcommand{\R}{\mathbb{R}}
\newcommand{\Sm}{\mathbb{S}}
\newcommand{\A}{\mathcal A}
\newcommand{\wu}{\widehat{u}}

\numberwithin{equation}{section}

\title[Trudinger-Moser type inequality in 2D]{New analytical and geometrical aspects on Trudinger-Moser type inequality in 2D}
\author{Natalino Borgia}
\author{Silvia Cingolani}
\author{Gabriele Mancini}

\begin{document}

\begin{abstract}
The present survey is devoted to results on Trudinger-Moser inequalities in two dimension. We give a brief overview of the history of these celebrated inequalities and, starting from the geometric problem that motivated Moser's original work, we discuss the connection between Onofri's inequality for the unit sphere and sharp inequalities on Euclidean domains. 
 Finally, we present recent results 
 and new insights into
nonlocal interaction energy functionals in two dimension, involving logarithmic kernels. 
\end{abstract}

\maketitle

\section{Introduction}

\bigskip	

If $\Omega$ is a bounded open subset of $\re^2$, the Sobolev space 
$H^1_0(\Omega)$ is a borderline case for the Sobolev embedding since it is
continuously embedded into $L^q(\Omega)$ for any
$1 \leq q < \infty$, but 
$
H^1_0(\Omega) \not \hookrightarrow L^\infty(\Omega).
$
If $\Omega$ is the unit disk, a counterexample is given by the function
$
u(x)= \ln(1- \ln|x|).
$

In 1967, Trudinger \cite{Tru} proved that $H^1_0(\Omega)$ is embedded into the Orlicz space 
$L_{\phi}(\Omega)$  corresponding to the Young measure $\phi(t)= e^{t^2} -1$ (see also Yudovich  \cite{Yudo} and Pohozaev \cite{poho}). In 1971, Moser \cite{M} established a  celebrated result which we recall for the two-dimensional case.

\begin{theorem}[\textbf{Trudinger-Moser inequality in 2D}]\label{MTRn}
	\smallskip
	Let $\Omega \subset \mathbb{R}^2$ be a bounded domain. Let $\mathcal{H}:=\{u \in H^1_0(\Omega) \ | \  \lVert \nabla u \rVert_2 \leq 1 \}$, where $\ \lVert\, \cdot\,\rVert_2$ denotes the usual $L^2$ norm. Then
	\begin{itemize}
		\item[$1)$] $\displaystyle \forall \, u \in H^1_0(\Omega)$ : $\quad e^{ u^{2}} \in L^1(\Omega). $
		\item[$2)$] There exists a constant $C \geq 1$ depending on the dimension only, such that 
		\begin{equation}\label{dis}
		\displaystyle \sup_{u \in \mathcal H}  \frac{1}{|\Omega|} \int_{\Omega} e^{4 \pi u^2 }\, dx  \leq  C.
		\end{equation}			 
		\small
		
		\item[$3)$] If $\alpha > 4 \pi,$ then $$ \displaystyle \sup_{u \in \mathcal H}  \frac{1}{|\Omega|} \int_{\Omega} e^{\alpha u^{2} }\, dx = + \infty.$$
	\end{itemize}
\end{theorem}

From Theorem \ref{MTRn}, it follows that
the maximal allowed growth is of exponential type and that
for any {$p>1$} there exists a sequence $(u_\epsilon)_{\epsilon>0}$ of functions in $H^1_0(\Omega)$ with $\|\nabla u_\epsilon \|_2=1$,
verifying {$\int_\Omega e^{4\pi p u_\epsilon^2 }\, dx \to +\infty$} as $\epsilon \to 0$. On the other hand, for all $u \in H^1_0(\Omega)$,   $\int_\Omega e^{4 \pi p u^2 }\, dx < \infty$ for all $p> 0$.

In 1985, Lions  \cite{Lions} proved that if $(u_\epsilon)$  is a sequence 
of functions in $H^1_0(\Omega)$ with $\|\nabla u_\epsilon \|_2=1$ such that $u_\epsilon \to u_0$ weakly in $H^1_0(\Omega)$, then for any $p < 1/ (1- \|\nabla u_0\|^2_2)$  it results $$\limsup_{\epsilon \to 0}  \int_\Omega e^{4 \pi p u_\epsilon^2 }\, dx  < + \infty.$$
Of course this fact gives extra information if   the sequence $u_\epsilon$ has a weak limit $u_0 \neq 0$.

\smallskip
A challenging problem on Trudinger-Moser inequalities was whether the extremal functions	exist or not. It has remained open for quite some years.
In 1986,
Carleson and Chang \cite{carlesonchang} proved that 
the supremum is attained in $(\ref{dis})$ if $\Omega = B_1$ is the unit disc in $\re^2$,
by using a symmetrization argument.
Successively
Struwe  \cite{Struwe} showed  that 
the supremum is attained in $(\ref{dis})$ if $\Omega$ is close to a ball in the sense of measure. Flucher \cite{flucher} extended this result to arbitrary
bounded domains in $\R^2$ using conformal deformations.

In 2004 Adimurthi and Druet \cite{AD}
highlighted that the result of Lions \cite{Lions} suggested 
\begin{equation}\label{disadi1}
	C_\alpha (\Omega):= 	\sup_{u \in H_0^1(\Omega), \|\nabla u\|_2\leq 1} \int_\Omega e^{4 \pi u^2(1 + \alpha \|u\|^2_2)}  dx < +
	\infty
\end{equation} 
for some $\alpha>0$.
In fact, they proved that if $\Omega$ is a smooth bounded domain in $\re^2$, denoted by 
$\lambda_1(\Omega)$  the first eigenvalue of the Laplacian with respect to the Dirichlet boundary condition, it results that
 for any 
$0 \leq \alpha < \lambda_1(\Omega)$, $C_\alpha(\Omega)< + \infty$. On the contrary, for any $\lambda \geq \lambda_1(\Omega)$, $C_\alpha(\Omega)= + \infty$.
These results in \cite{AD} add new information  when the weak limit of $u_\epsilon$ is $0$. 
A stronger version  of $(\ref{disadi1})$ was obtained by Tintarev \cite{Tintarev} and  existence of extremal functions was proved  by Yang \cite{Yang}. { We point out that the existence of extremal functions for these improved inequalities is an extremely delicate problem. For example,  the  supremum in \eqref{disadi1} is attained  for small values of $\alpha$, but it is not attained when $\alpha$ is close to $\lambda_1(\Omega)$, see \cite{ManThi}. In general, even small perturbations of the exponential term in \eqref{dis} may affect the attainability of the supremum, as it is pointed out by \cite{pruss} and clarified by \cite{MM} and \cite{IMNS}. We refer to \cite{Thizy} for recent developments on the subject.  }

As regards Trudinger-Moser inequality $(\ref{dis})$ on unbounded domains, we mention 
for the subcritical growth  $\alpha < 4 \pi$ the results by Cao \cite{cao} and Adachi, Tanaka \cite{AT} (see also \cite{doO}).
In  the critical case $\alpha =4 \pi$,  the inequality with the sharp exponent (proved by Moser) fails on the whole plane. Ruf (see \cite{ruf} and \cite{LiRuf}) proved that 
 there exists a constant $d>0$  such that, for any $\Omega$ unbounded domain in $\re^2$
\begin{equation}\label{dis1}
	\sup_{u \in H_0^1(\Omega), \|u\|\leq 1} \int_\Omega \bigl(e^{4 \pi u^2} -1 \bigr) dx \leq d,
\end{equation}
where  $\|u\|= \bigr(\int_\Omega (|\nabla u|^2 + u^2)dx \bigl)^{1/2}$ denotes the usual Sobolev norm on $H_0^1(\Omega)$.
 
The supremum  is attained if $\Omega=\R^2$ {(see \cite{LiRuf,Ishi}) and  the inequality is sharp: for any growth $e^{\alpha u^2}$} with $\alpha > 4 \pi$, the corresponding supremum equals $+ \infty$. We mention \cite{LiRuf} for a general dimension and \cite{IMN} for
a more recent improvement of the inequality on the whole space $\re^2$. Furthermore an equivariant version  of (\ref{dis1}) has been obtained in \cite{CST} and  sharp weighted Trudinger-Moser type inequalities were proved in \cite{DongLamLu,CaRuf}.
{For unbounded domains with  $\lambda_1(\Omega)>0$, bounds on the full Sobolev norm are unnecessary, as one can prove that
$$
\sup_{u\in \mathcal H} \int_{\Omega} (e^{4\pi u^2}-1)\,dx <\infty,
$$  
where $\mathcal H$ is defined as in Theorem \ref{MTRn}, see  \cite{ManSand} and \cite{BatMan}.}

{Even if in this survey we focus on inequalities for $H^1_0(\Omega)$, it is important to mention the work of Adams \cite{Adams}, who proved similar inequalities for  domains $\Omega\subseteq \R^n$ and the Sobolev spaces $W^{k,p}(\Omega)$ whenever $k p= n$. We also refer to \cite{Martinazzi} for Trudinger-Moser-Adams type inequalities on fractional order Sobolev spaces. Concerning the study of extremal functions for this more general family of inequalities, we refer to \cite{YoungAdams,ManDLT,MM2} and references therein.}

{Trudinger-Moser type inequalities have also been widely  studied  for Sobolev spaces on Riemannian surfaces and manifolds. Indeed, as we discuss in the next section, they have important applications to prescribed curvature problems in conformal geometry. A general formulation of Trudinger-Moser-Adams inequalities for compact Riemannian manifolds was given by Fontana in \cite{Fontana} ( see also  \cite{Li1,Li2} for related results). } Recently subcritical and critical Trudinger-Moser inequalities are proved for complete noncompact Riemannian manifolds in \cite{LiLu}.

The singular Trudinger–Moser inequality and its extremal functions were considered in \cite{ChenSing,Troy,ASandeep,CR,IM,YangZhu} and generalizations to the whole space in
\cite{Yang,LiYang}.

\smallskip
The aim of the present survey is two-fold. On the one hand we discuss  some connections between different notions of  Trudinger-Moser inequality in 2D, announcing some forthcoming results in \cite{BCM}. In particular, in section \ref{SecOnofri}, starting from geometrical and physical applications, we introduce an inequality due to Onofri  \cite{O} (see \eqref{Onofri})  for Sobolev spaces on the unit sphere and we discuss its connection with similar inequalities for domains in $\R^2$.

 On the other hand  we present some recent results and new insights into
 nonlocal interaction energy functionals in two dimension, arising from  mathematical
 models for chemotaxis \cite{Dolbeault-Perthame}, in the statistical mechanics of selfgravitating clouds \cite{Suzuki,W} 
 and in the description of vortices in turbulent Euler flows \cite{CLMP}.
 In particular, in section 3, we state some 
 Trudinger-Moser type inequalities with logarithmic kernel in \cite{CiWe2} and announce some forthcoming developments in \cite{CiWeYu}.

\section{Onofri-type inequalities and geometric connections}\label{SecOnofri}

The Trudinger-Moser inequality has applications to several mathematical fields ranging from mathematical physics to differential geometry and calculus of variations. 
Among the reasons that led Moser to improve Trudinger's result, it is important to mention  Nirenberg's problem on $\Sm^2$, which consists in describing the set of smooth functions that can be obtained as Gaussian curvature of a Riemannian metric conformally equivalent to the standard spherical metric. 

 Let $\Sigma$ be a smooth surface, we recall that two Riemannian  metrics $g_0$ and $g$ on $\Sigma$ are said to be pointwise conformally equivalent if there exists  $u\in C^\infty(\Sigma)$ such that $g = e^u g_0$ on $\Sigma$. In this setting, the following identity holds:
$$
\Delta_{g_0} u + 2K_g e^{u} = 2 K_{g_0},
$$ 
where $\Delta_{g_0}$ denotes the Laplace-Beltrami operator on $\Sigma$, and $K_g$ and $K_{g_0}$ are the Gaussian curvatures of $\Sigma$ corresponding respectively to $g$ and $g_0$. When $\Sigma = \Sm^2$ and $g_0$ is the standard spherical metric, then a function $K \in C^\infty(\Sm^2)$ is the Gaussian curvature of a metric which is pointwise conformally equivalent to $g_0$ if and only if there exists $u\in C^\infty(\Sm^2)$ satisfying 
\begin{equation}\label{EqS2}   \Delta_{g_0}u + 2K e^u=2 \quad \mbox{ on } \quad (\Sm^2,g_0).
\end{equation}
It is possible to give a variational formulation of \eqref{EqS2} by considering the functional 
$$ \displaystyle \Phi(u) = \frac{1}{2} \int_{\Sm^2} |\nabla_{g_0} u|^2_{g_0} \, dv_{g_0} + 2\int_{\Sm^2} u  \, dv_{g_0}  - 8 \pi \ln \left( \frac{1}{4\pi}\int_{\Sm^2} K e^u\, dv_{g_0} \right),$$ 
which is well defined in $$ \A=: \left\{ u \in H^1(\Sm^2,g_0) \quad | \quad \int_{\Sm^2} K e^u \, dv_{g_0} >0   \right\},$$
provided ${\displaystyle \max_{\Sm^2} K >0}$. A function $u \in H^1(\Sm^2,g_0)$ is a weak solution to \eqref{EqS2} if and only if it is a critical point of $\Phi$ and 
$$
\int_{\Sm^2} K e^{u}dv_{g_0}  = 4\pi.
$$
As an analog of Theorem \ref{MTRn} for the Sobolev space $H^1(\Sm^2,g_0)$, Moser proves that there exists a constant $S>1$ such that for any $u \in H^1(\Sm^2,g_0)$ with $\|\nabla_{g_0} u\|_{L^2(\Sm^2,g_0)}\le 1$:
\begin{equation}\label{MTS2}\displaystyle  \frac{1}{4 \pi} \int_{\Sm^2} e^{4 \pi {(u-\overline{u})^2}} dv_{g_0} \leq S, \quad \text{where} \quad \overline{u} := \frac{1}{4\pi}\int_{\Sm^2} u \, dv_{g_0}. 
\end{equation}

As in the Euclidean inequality 	\eqref{dis}, the exponent $4\pi$ in \eqref{MTS2} is sharp. As a straightforward consequence of \eqref{MTS2} and Young's inequality, one can prove the weaker  inequality:
\begin{equation}\label{ConsMS2}
\ln\left(\frac{1}{4\pi}\int_{\Sm^2} e^u \, dv_{g_0} \right) \leq  \frac{1}{16 \pi} \int_{\Sm^2} | \nabla_{g_0} u |_{g_0}^2 \, dv_{g_0} + \frac{1}{4 \pi} \int_{\Sm^2}u \, dv_{g_0} + { \ln S},
\end{equation}
for any $u \in H^1(\Sm^2,g_0),$ where $S$ is the constant appearing in \eqref{MTS2}, and the coefficient $\frac{1}{16\pi}$ is sharp. In particular,  \eqref{ConsMS2} shows that $\Phi$ is bounded from below in $\A$. While this is not sufficient to apply standard minimization techniques (actually, it can be shown that minimum points for $\Phi$ in $\A$ do not exist unless $K$ is constant), the existence of critical points for $\Phi$ can be obtained via improved versions of \eqref{ConsMS2} under suitable assumptions on $K$. For example, in \cite{M2} Moser proves that if $u\in H^1(\Sm^2,g_0)$ is even, then \eqref{ConsMS2} can be replaced by the stronger inequality:  
\begin{equation}\label{ConsMEven}
\ln \left(\frac{1}{4\pi}\int_{\Sm^2} e^u \, d\nu_{g_0} \right) \leq  \frac{1}{32 \pi} \int_{\Sm^2} | \nabla_{g_0} u |_{g_0}^2 \, d\nu_{g_0} + \frac{1}{4 \pi} \int_{\Sm^2}u \, d\nu_{g_0} + C ,
\end{equation}
for some $C\in \R$. In particular, if $K\in C^\infty(\Sm^2)$ is even, the restriction of $\Phi$ to the space of even functions in $H^1(\Sm^2,g_0)$ with zero average is coercive, and one can find a critical point of $\Phi$ by means of Palais symmetric criticality principle. This lead to one  of first general existence results for the Nirenberg problem: if $K\in C^\infty(\Sm^2)$ is an even function and $\displaystyle{\max_{\Sm^2} K>0}$, then there exists a solution of \eqref{EqS2}. Over the years, many improved versions of \eqref{ConsMS2} have been discovered. For example, we mention that \eqref{ConsMS2} can be improved for functions with vanishing barycenter (see \cite{Aubin} and \cite{MalRuiz}) or for functions such that mass of the measure $e^{u}dv_{g_0}$ is spread in at least two disjoint open subsets of $\Sm^2$ (see \cite{CL,Malchiodi}). We refer to {\cite{Berger,KW,CY1,CY2,StruweCurv} and references therein for further results concerning Nirenberg's problem on $\Sm^2$.}

Another interesting aspect of \eqref{ConsMS2} is that while the coefficient  $\frac{1}{16\pi}$ is sharp, the constant $\ln S$ is not. Indeed, Onofri in \cite{O} proved that
$$
\inf_{H^1(\Sm^2,g_0)} \Phi  = 0 
$$
and he classified all minimum points of $\Phi$. Namely, we have 
\begin{equation}\label{Onofri}
\ln \left(\frac{1}{4\pi}\int_{\Sm^2} e^u \, dv_{g_0} \right) \leq  \frac{1}{16 \pi} \int_{\Sm^2} | \nabla_{g_0} u |_{g_0}^2 \, dv_{g_0} + \frac{1}{4 \pi} \int_{\Sm^2}u \, dv_{g_0}, 
\end{equation}
for any $\, u \in H^1(\Sm^2,g_0)$, and equality  holds if and only if $u  = \ln |J_\psi| + C$, where $\psi:\Sm^2 \rightarrow \Sm^2$ is a conformal diffeomorphism of $(\Sm^2,g_0)$, $J_\psi$ denotes the Jacobian of $\psi$, and $C\in \R$. In \cite{DEJ}, \eqref{Onofri} is obtained as limit of Sobolev's inequalities for $W^{1,p}(\R^2)$ as $p\nearrow 2$.  Onofri's inequality \eqref{Onofri} plays a role in spectral analysis for the Laplace-Beltrami operator. Indeed, thanks  to Polyakov's formula (see \cite{Poly1,Poly2,OPS,OPS2}), we have 
$$
\frac{\det (-\Delta_g)}{\det(-\Delta_{g_0})} =e^{-\frac{1}{24}\Phi(u)},
$$
for any metric $g = e^u g_0$ with $u \in C^\infty(\Sm^2)$. In particular,  Onofri's inequality \eqref{Onofri} implies that, in the conformal class of $g_0$, the determinant of the Laplace-Beltrami operator is maximized by  $g_0$.

An analog of \eqref{Onofri} for the unit disk in $\R^2$ can be deduced from the work of Carleson and Chang \cite{carlesonchang}. As a crucial step in their proof of the existence of extremal functions for \eqref{dis}, they prove an inequality which is equivalent to:  
\begin{equation}\label{OnofriDisc}
\ln \left(\frac{1}{\pi}\int_{B_1} e^u \,dx \right) <  \frac{1}{16\pi}\int_{B_1} |\nabla u|^2 dx  + 1,
\end{equation}
for any $u\in H^1_0(B_1)$. Again, this inequality is sharp in the sense that
$$
\inf_{u\in H^1_0(B_1)} \left( \frac{1}{16\pi}\int_{B_1} |\nabla u|^2 dx - \ln \left(\frac{1}{\pi}\int_{B_1} e^u \,dx \right) \right) = -1,
$$
and 
$$
\inf_{u\in H^1_0(B_1)}\left( \alpha \int_{B_1} |\nabla u|^2 dx - \ln \left(\frac{1}{\pi}\int_{B_1} e^u \,dx \right) \right) = -\infty,
$$
for any $\alpha < \frac{1}{16\pi}$. More generally, it was proved in \cite{CLMP} (see also \cite{CCL}) that if $\Omega\subseteq \R^2$ is a bounded domain with smooth boundary,  then for any $u \in H^1_0(\Omega)$  
\begin{equation}\label{OnofriDomain}
\ln \left(\frac{1}{|\Omega|}\int_{\Omega} e^u \,dx \right) \le  \frac{1}{16\pi}\int_{\Omega} |\nabla u|^2 dx  + 1 + {4\pi \sup_{\Omega}\gamma(x) + \ln \frac{\pi}{|\Omega|}},
\end{equation}
where $\gamma$ is Robin's function of $\Omega$, which is defined for  $x\in \Omega$ as $\gamma (x) = H_x(x)$, where $H_x$ is the unique solution of 
$$
\begin{cases}
\Delta H_x = 0  & \text{ in }\Omega,\\
H_x(y) = \frac{1}{2\pi }\ln|x-y|  & \text{ on }\partial \Omega. 
\end{cases}
$$
Differently from \eqref{OnofriDisc},  \eqref{OnofriDomain} is  not always sharp. Indeed, denoting 
$$
I(\Omega):= \inf_{u\in H^1_0(\Omega)} \left( \frac{1}{16\pi}\int_{\Omega} |\nabla u|^2 dx - \ln \left(\frac{1}{|\Omega|}\int_{\Omega} e^u \,dx \right)\right),
$$
there are domains $\Omega$ (see \cite[Section 7]{CLMP0}) for which 
\begin{equation}\label{nonsharp}
 I(\Omega)>- 1- {4\pi \sup_{\Omega}\gamma(x) - \ln \frac{\pi}{|\Omega|}}.
\end{equation}
The domains for which \eqref{OnofriDomain} is sharp are known in the literature as domains of first kind, while the ones for which  \eqref{nonsharp} holds are known as domains of second kind. Inequality \eqref{OnofriDomain} plays a crucial role in the study of the mean field equation 
\begin{equation}\label{mean}
-\Delta u  = \frac{\rho e^u}{\int_{\Omega}e^u \,dx} \quad \text{ in }\Omega,
\end{equation}
which appears in the statistical mechanic description of vortex formations in 2d-models for turbulent flows in fluid dynamics (see e.g. \cite{JM,CLMP}). We refer to \cite{CL1,CL2,Malchiodi,Bart,BarMal} for existence results  for  \eqref{mean} and related problems.

\medskip
In \cite{M} Moser explicitly states that it seems impossible to  deduce \eqref{dis} from \eqref{MTS2} or \eqref{MTS2} from \eqref{dis}. The proofs of the two inequalities are somehow independent, although they are based on similar arguments relying on symmetrization techniques. 
%
%
 It turns out that a more explicit connection exists between \eqref{Onofri} and \eqref{OnofriDisc}. Indeed, one can show that  \eqref{Onofri} and \eqref{OnofriDisc} are  equivalent by exploiting the conformal equivalence between  $\R^2$ and $\Sm^2$ minus one point. In the remaining part of this section, we sketch a simple argument to deduce \eqref{OnofriDisc} from \eqref{Onofri}, which is inspired by \cite{IM}, where the authors consider generalized versions of \eqref{OnofriDisc} and \eqref{Onofri} involving singular weights. {In a forthcoming paper \cite{BCM}, we will give a complete proof of the equivalence between \eqref{Onofri} and \eqref{OnofriDisc} and investigate extensions to  dimension $n\ge 3$ (see \cite{DD}, \cite{ABG} and \cite{LamLu}). }
 
 \medskip
 In the following, we assume that \eqref{Onofri} holds in $H^1(\Sm^2,g_0)$. We take a function  $u\in H^1_0(B_1)$, and we prove that \eqref{Onofri} holds.  We denote by $N=(0,0,1)$ the north pole of $\Sm^2$ and we let $\pi:\Sm^2 \setminus\{N\} \rightarrow \R^2$ be the standard stereographic projection. We shall consider the function $$v_0(x) = 2\ln \left(\frac{2}{1+|x|^2}\right).$$
We stress that $v_0$ satisfies 
 \begin{equation}\label{eqv0}
- \Delta v_0 = 2 e^{v_0} \quad \text{ on } \R^2. 
\end{equation}
Moreover, for any $f:\Sm^2 \rightarrow \R$, $f\in L^1(\Sm^2,g_0),$ we have that 
$$
\int_{\Sm^2} f \,dv_{g_0} = \int_{\R^2} (f\circ \pi^{-1})\, e^{v_0} \,dx.
$$
For any $r>0$, we define $u_r: \mathbb{R}^2 \to \mathbb{R}$ as 
$$ \displaystyle 
u_r(x)=
\begin{cases}
\displaystyle u \left( \frac{x}{r} \right) - v_0(x) & \text{if $|x| \leq r$,} \bigskip \\
-2\ln\left(\frac{2}{1+r^2}\right) & \text{if $|x| > r$,}
\end{cases}
$$
and  $\wu_r: \Sm^2  \to \mathbb{R}$ as 
$$ \displaystyle 
\wu_r(x)=
\begin{cases}
u_r \circ \pi(x)  & \text{ if } x\neq N\\
-2\ln\left(\frac{2}{1+r^2}\right) &  \text{ if } x= N.
\end{cases}
$$
Observe that ${\wu_r}$ is constant in a neighborhood  of $N$ and ${\wu_r} \in W^{1,2}(\mathbb{S}^2,g_0)$. Then, we can apply \eqref{Onofri} to get:
$$ \displaystyle \ln \left( \frac{1}{4 \pi} \int_{\mathbb{S}^2} e^{{\wu_r}} \, dv_{g_0} \right) \leq  \frac{1}{16 \pi} \int_{\mathbb{S}^2} | \nabla_{g_0} {\wu_r} |_{g_0}^2 \, dv_{g_0} + \frac{1}{4 \pi} \int_{\mathbb{S}^2} {\wu_r} \, dv_{g_0}.
$$
A direct computation shows that 
\begin{align*}
\displaystyle 
\frac{1}{4 \pi} \int_{\mathbb{S}^2} e^{\widehat{u_r}} \, dv_{g_0} = \frac{r^2}{4 \pi} \int_{B_1} e^{u} \, dx  + \frac{1+r^2}{4}.
\end{align*}
Similarly, we have
\begin{equation}\label{media}
\displaystyle  \frac{1}{4 \pi} \int_{\mathbb{S}^2} {\wu_r} \, dv_{g_0}= \frac{1}{4 \pi} \int_{B_r} u_r e^{v_0} \, dx  + 2 - \ln 4 + o_r(1),
\end{equation}
and 
\begin{equation}\label{gradv0} \displaystyle  \frac{1}{16 \pi} \int_{B_r} |\nabla v_0|^2 \, dx = 2\ln r - 1+ o_r(1),
\end{equation}
where $o_r(1)\to 0$ as $r \to +\infty$. 
Using \eqref{eqv0} and \eqref{gradv0} and \eqref{media}, we find
$$ \displaystyle 
\begin{aligned}\frac{1}{16 \pi} &\int_{\mathbb{S}^2} | \nabla_{g_0} {\wu_r} |_{g_0}^2 \, dv_{g_0} \\& =\frac{1}{16 \pi} \int_{B_r} \left| \nabla u_r \right|^2 \, dx -\frac{1}{4 \pi} \int_{B_r} u_r e^{2v_0} dx + 2\ln r - 1+o_r(1) \\
&= \frac{1}{16 \pi} \int_{D} \left| \nabla u \right|^2 \, dx - \frac{1}{4 \pi} \int_{\mathbb{S}^2} {\wu_r} \, dv_{g_0} + \ln\left( \frac{r^2}{4}\right) + 1+o_r(1).
\end{aligned}
$$
Summing up, we have shown that
$$ \displaystyle \ln \left(\frac{1}{\pi} \int_{B_1} e^{u} \, dx + \frac{1+r^2}{r^2} \right) \leq \frac{1}{16 \pi} \int_{D} \left| \nabla u \right|^2 \, dx +
 1  + o(1).$$
Passing to limit as $r \to + \infty,$ we conclude
$$ \displaystyle \ln \left(\frac{1}{\pi} \int_{B_1} e^{u} \, dx \right) {<} \ln \left(\frac{1}{\pi} \int_{B_1} e^{u} \, dx + 1 \right) \leq \frac{1}{16 \pi} \int_{B_1} \left| \nabla u \right|^2 \, dx  +1.$$
Hence $u$ satisfies \eqref{OnofriDisc}.


\section{Trudinger-Moser type inequalities with logarithmic kernels}

In this section we consider some nonlocal interaction functionals, which naturally
arise in  2D mathematical models for chemotaxis and in other physical contexts, as statistical dynamics, vortex theory.

These two-dimensional problems have remained for a long time a quite open field of study, due to the difficulties to handle changing sign kernels,  not bounded from above and below. In literature the first results are just numerical.  Recently variational results have been employed to treat integro-differential equations, involving  logarithmic convolution potentials,
where classical PDE theory fails because of the presence of the non-locality
(see 
\cite{BoCiSe,BoCiVa,BVS,CiJe,CiWe,DuWe,Masaki,Masaki2,stubbe,SV}).

In particular
we consider the following class of free energy functionals  
$$u \mapsto \int_{\Omega} \int_{\Omega} \ln \frac{1}{|x-y|}F(u(x))F(u(y))\,dx dy,$$
where 
$\Omega=B_1$ is the unit ball of $\re^2$
 and 
$F: \R \to \R$ grows exponentially.

In \cite{CiWe2} we faced the following challenging questions:

\begin{itemize}
	\item{\sl Moser-Trudinger type inequalities in presence of a logarithmic convolution potential;}
		
	\item{\sl characterization of the critical nonlinear growth rates for these inequalities;}
	
	\item{\sl associated Euler-Lagrange equation;}
	
	\item{\sl qualitative properties of the maximizers;}
	
	\item{\sl symmetry and uniqueness results.}
\end{itemize}


Here we give some outlines of the results and the proofs 
in \cite{CiWe2}. Firstly we consider the problem of maximizing 
the quantity 
	\begin{equation}
		\label{eq:logarithmic-convolution-kernel}
		\Phi(u) := \int_{B_1} \int_{B_1} \ln \frac{1}{|x-y|}F(u(x))F(u(y))\,dx dy 
\end{equation}
among functions
in the set 
{$$
\cB:= \{u \in H^1_0(B_1)\::\: |\nabla u|_2 \le 1 \}.
$$ }

We make the following general assumption on $F$:
\smallskip 

\begin{itemize}
	\item[$(A_0)$] 
	$F: \R \to [0,\infty)$ is even, continuous, and strictly increasing on $[0,\infty)$. Moreover, there exist constants $\alpha,c>0$ with 
		\begin{equation}
			\label{eq:general-growth-condition}
			F(t) \le c e^{\alpha t^2}  \qquad \text{for $t \in \re$.}
	\end{equation}
\end{itemize}

Since the logarithmic kernel changes sign, it is not a priori clear if the double integral has a well-defined value. 
To clarify this point,	we split the kernel $\ln \frac{1}{|\,\cdot\,|}$ into its positive and negative part and define the functionals 
$$
\begin{aligned}
		\Phi_+(u)&=  \int_{B_1} \int_{B_1} \ln^+ \!\frac{1}{|x-y|}\, F(u(x))F(u(y))\,dx dy \in [0, +\infty], \\ 
		\Phi_-(u)&=  \int_{B_1} \int_{B_1} \ln^+\! |x-y|\, F(u(x))F(u(y))\,dx dy \in [0, + \infty],
\end{aligned}
$$
where $u \in H^1_0(B_1)$, and $\ln^+ = \max \{\ln,0\}$. 
By the classical Trudinger-Moser inequality,  	
\begin{equation*}
	\label{general-nonuniform-growth-estimate}
		\int_{B_1}e^{\alpha u^2}dx < \infty  \qquad \text{for any $u \in H^1_0(B_1)$ and any $\alpha>0$.}
\end{equation*} 
By the growth condition
$$F(t) \le c e^{\alpha t^2}  \qquad \forall t \in \R,$$ we infer 
$F(u) \in L^s(\re^2)$ for any 
$u \in H^1_0(B_1)$, $s \in [1,\infty)$.
Since 
\begin{equation*}
	\ln^+ \frac{1}{|\cdot|} \in L^s(\re^2) \qquad \text{for $s \in [1,\infty)$,}  
\end{equation*}
choosing $s=2$, we deduce by Young convolution inequality that $\ln^+\frac{1}{|\cdot|} * F(u) \in L^\infty(\re^2)$ and
thus $\Phi_+(u)<\infty$, for every $u \in H^1_0(B_1)$.

\smallskip
Moreover, we have
$$
[\ln^+|\cdot| * F(u)](x) = \int_{B_1}\ln^+|x-y|F(u(y))\,dy < \ln 2 |F(u)|_1 \qquad \text{ $\forall x \in B_1$}
$$
and therefore $\ln^+|\cdot| * F(u) \in L^\infty(B_1)$ and thus $\Phi_-(u)<\infty$ for every $u \in H^1_0(B_1)$.

Taking into account that 
$\displaystyle \ln \frac{1}{|\cdot|} =\ln^+ \frac{1}{|\cdot|} -	\ln^+ {|\cdot|},$
we have 
	\begin{equation*}
		\Phi(u):=\int_{B_1} \int_{B_1} \ln \!\frac{1}{|x-y|}\, F(u(x))F(u(y))\,dx dy = \Phi_+(u) -\Phi_-(u)  \in ]-\infty,\infty[ 
\end{equation*}
for every $u \in H^1_0(B_1)$.

\medskip
Now in order to study  the maximization problem related to the quantity
\begin{equation*}
		\label{eq:def-Phi-new}
		\Phi(u):=\int_{B_1} \int_{B_1} \ln \!\frac{1}{|x-y|}\, F(u(x))F(u(y))\,dx dy 
\end{equation*}
within the set 
$$
\cB:= \{u \in H^1_0(B_1)\::\: |\nabla u|_2 \le 1 \}.
$$ 
we need to distinguish different forms of asymptotic growth of the nonlinearity $F$. As in Trudinger-Moser inequality, the value $4\pi$ will play a key role.

\smallskip
Firstly we notice that if $F$ has a subcritical exponential growth $\alpha < 4 \pi$, 
an immediate proof of the claim 
	$$
m_1(F) := \sup_{u \in \cB} \Phi(u)   <\infty
$$
can be given using  the logarithmic Hardy-Littlewood-Sobolev inequality 
$$
\int_{B_1} \int_{B_1} \ln \frac{1}{|x-y|}f(x)f(y)\,dx dy 
	\le \frac{|f|_{1}}{2}\Bigl(|f|_1 \bigl(c_0 +\bigl|\ln |f|_1\bigr|\bigr) + \int_{B_1} f \ln f\,dx\Bigr),
$$
which holds for any a.e. non-negative function $f\in L^1(B_1)$ with $f \ln f\in L^1(B_1)$, see \cite[Theorem 2]{Bek}.

Indeed,  let
$u \in \cB \cap L^\infty(B_1)$, and let 
$v:= F(u)$.
Then $v \le c_1 e^{\alpha u^2}$, $\alpha< 4\pi$, $c_1>0$
and 
$$
v \ln v \le c_1 e^{\alpha u^2} \Bigl(\alpha u^2 + \ln c_1 \Bigr) \le c_2 e^{4 \pi u^2}
$$
with $c_2>0$. By the classical Trudinger-Moser inequality,
we have 
\begin{equation*}
	\label{eq:remark-intro-3}
	|v|_1 \le c_1 \int_{B_1}e^{4 \pi u^2}dx \le c_1 c(B_1) =:c_2, \quad \qquad\int_{B_1} v \ln v \,dx \le c_2 c(B_1)=: c_3. 
\end{equation*}
Combining with the 
logarithmic Hardy-Littlewood-Sobolev inequality	we infer
\begin{align*}
	&	\Phi(u)= \int_{B_1} \int_{B_1} \ln \frac{1}{|x-y|}v(x)v(y)\,dx dy 
	\\ &\le \frac{|v|_{1}}{2}\Bigl(|v|_1 \bigl(c_0 +\bigl|\ln |v|_1\bigr|\bigr) + \int_{{B_1}} v \ln v\,dx\Bigr)
	\le \frac{c_2}{2}\bigl(c_2 \bigl({c_0} + |\ln c_2|\bigr) + c_3\bigr)<\infty.
\end{align*}	
For general $u \in \cB$, the same inequality follows by approximation and  $m_1(F) < \infty$.
This argument does not apply if $F$ is critical, so it does not allow for sharp critical growth rates.

\smallskip

In \cite{CiWe2} we define some  critical rates for the Log-TM inequality.

\begin{definition}
	Let 
	$\beta \in \re$, and let $g: \re \to  \re$  be an arbitrary function. We say that $g$ has 
	\smallskip
	\begin{itemize}
		\item[(i)] at most $\beta$-critical growth if $|g(s)| \le c\ e^{4\pi s^2}(1+|s|)^\beta$ for $s \in \re$ with some constant $c>0$;
		
		\smallskip
		\item[(ii)] at least $\beta$-critical growth if there exist $s_0, c>0$ with the property that 
		$$
			|g(s)| \ge c\ e^{4\pi s^2}|s|^\beta \qquad \text{for $|s| \ge s_0$.}
		$$        
	\end{itemize}
\end{definition}

In  \cite{CiWe2} we obtained the following result.

\bigskip
\begin{theorem} 
	\label{sec:introduction-main-thm-C-F2}
	Suppose that $F$ satisfies $(A_0)$.
	
	\begin{enumerate}
		\item If $F$ has at most $\beta$-critical growth
		for some 	$\beta \le -1$, then 
		$$
		m_1(F) := \sup_{u \in \cB} \Phi(u)   <\infty.
		$$
		\item If $F$ has at most 
$\beta$-critical growth  for some $\beta < -1$, then $m_1(F)$ is attained, and every maximizer for $\Phi$ in $\cB$ is, up to sign, a radial and radially decreasing function in $B_1$.  
		\item If $F$ has at least $\beta$-critical growth for some $\beta>-1$, then  $m_1(F)= \infty$.
	\end{enumerate}
\end{theorem}

\smallskip

	\medskip
	In order to maximize $\Phi$
	among the functions in $\cB:= \{u \in H^1_0(B_1)\::\: |\nabla u|_2 \le 1\}$,
	we consider
	$$
			\cB^*:=  \{u^* \::\: u \in \cB\}
	$$
	the corresponding  Schwartz symmetrized set. 
	Here $u^*$ is nonnegative, radial and nonincreasing in the radial variable.

	\smallskip  
	By the P\'olya-Szeg\"o inequality,
	for any $u \in H^1(\R^2)$ we have $u^* \in H^1(\R^2)$ and 
	$$
	|u^*|_p = |u|_p, \quad |\nabla u^*|_2 \le |\nabla u|_2,  \qquad \qquad \text{ $ \forall  p \in [2, +\infty[$.}
	$$
	Consequently, we have $u^* \in \cB$ if $u \in \cB$, and so
	$$\cB^* \subset \cB.$$
	
	Now  
	by Riesz's rearrangement inequality, 	for every $u \in \cB$
	\begin{align*}	
		\Phi_+(u) &=
		\int_{B_1} \int_{B_1} \ln^+ \!\frac{1}{|x-y|}\, F(u(x))F(u(y))\,dx dy  \\  & \leq 
		\int_{B_1} \int_{B_1} \ln^+ \!\frac{1}{|x-y|}\, [F(u(x))]^*[F(u(y))]^* \,dx dy \\ 
	& = 	\int_{B_1} \int_{B_1} \ln^+ \!\frac{1}{|x-y|}\, [F(u^*(x))] [F(u^*(y))] \,dx dy = \Phi_+(u^*)
	\end{align*}
	Therefore it is immediate that  $$m_1(F): = \sup_{\cB}\Phi(u) \le 
	m_1^+(F) : = \sup_{\cB}\Phi_+ = \sup_{u \in \cB^*} \Phi_+(u).$$

	In \cite{CiWe2} we prove that 
	$$\sup_{u \in \cB^*} \Phi_+(u)  < + \infty$$ and thus the	
	the following proposition holds.

	\begin{proposition}
		\label{C-finiteness}
		Suppose that $F$ satisfies $(A_0)$ and has at most $\beta$-critical growth for some ${\beta\le -1}$. Then we have 
		$$			m_1(F): =  \sup_{u \in \cB}\Phi(u) < \infty		$$
		
	\end{proposition}

\smallskip	
	We emphasize the optimality of the growth exponent $\beta = -1$ for the finiteness of $m_1(F)$. Indeed  we can define 
the Moser type functions $u_n \in H^1_0(B_1) \cap L^\infty(B_1)$, for any $n \in \na$, $n \geq 2$   
\begin{equation}\label{MF1}
u_n(x)= \frac{1} {\sqrt{2 \pi}} \frac{\ln(1/|x|)}{(\ln n)^{1/2}} \bigl(1 - 
\frac{1}{4 \ln n} \bigr)^{1/2}, \quad   \ \frac{1}{n} \leq |x| \leq 1,
\end{equation}
\begin{equation}\label{MF2}
u_n(x)= \frac{1} {\sqrt{2 \pi}} {(\ln n)^{1/2}} \bigl(1 - 
\frac{1}{4 \ln n} \bigr)^{1/2}, \quad  \  0 \leq |x| \leq \frac{1}{n}.
\end{equation}

A standard computation gives the following blow-up estimates. 
	
	\begin{proposition}
		\label{C-infiniteness}
		Suppose that $F$ satisfies $(A_0)$ and has at least $\beta$-critical growth for some 
				$\beta>-1$. Then, if $u_n$ defined as in \eqref{MF1}-\eqref{MF2},  we have $u_n \in \cB \cap L^\infty(B_1)$ and   	$\Phi(u_n) \to \infty$ as $n \to \infty$.
	\end{proposition}
	
Finally, using Concentration and Compactness Principle by Lions \cite{Lions}, we can recognize the existence of estremals for the Log-TM maximization problem when the parameter $\beta < -1$.

	\begin{proposition}
		\label{C-existence-of-maximizer}
		Suppose that $F$ satisfies $(A_0)$ and has at most $\beta$-critical growth with $\beta<-1$. Then the value 
		$m_1(F)< \infty$ is attained by a function $u \in \cB^*$. 
	\end{proposition}
	Let $(u_n)_n$ be a maximizing sequence in $\cB$ for $m_1(F)$. We may assume that $u_n \in \cB^*$ for $n \in \na$. Since $\cB$ is bounded in $H^1_0(B_1)$, we may also assume that  
	\begin{equation*}
		\label{eq:b-1-weak-convergence}
		u_n \weak u \in H^1_0(B_1) \qquad \text{with $u \in \cB^*$.}
	\end{equation*}
	\begin{proof}
	For reader's convenience we give a sketch of the proof in \cite{CiWe2}.
	By Lions \cite{Lions}, we have two possibilities. Either 
	\smallskip
	\begin{itemize}
		\item[i)] $u=0$ or 
		\item[ii)] $u\not = 0$, and $\int_{B_1}  e^{(4 \pi +t )u_n^2 } \ dx$ is bounded for some $t>0$ and thus 
		\begin{equation*}
			\label{eq:ii-l-1-convergence}
			e^{4 \pi u_n^2} \to e^{4 \pi u^2 }  \quad \hbox{in}  \ \ L^1(B_1).
		\end{equation*}
	\end{itemize}
	
	Assume first that $i)$ holds. In this case, we can prove   that  
	\begin{equation*}
	m_1(F) = \lim_{n \to \infty} \Phi(u_n) = \Phi(0).
	\end{equation*}
	On the other hand, 
	if $F$ satisfies $(A_0)$, it follows from that 
		\begin{align*}
			\Phi(u) &= 2 (2\pi)^2 \int_{0}^1 r F(u(r)) \ln \frac{1}{r} \int_{0}^r \rho  F(u(\rho)) d\rho dr \\ &> 2 (2\pi)^2 \int_{0}^1 r F(u(0)) \ln \frac{1}{r} \int_{0}^r \rho  F(u(0)) d\rho dr 
			=\Phi(0)
	\end{align*}	
	for every $u \in H^1_0(B_1) \setminus \{0\}$. Hence
	$$m_1(F) \in \bigl(\Phi(0), \infty \bigr]$$
	and therefore $m_1(F)$ is not attained at $u=0$.	So this case does not occur.

	It remains to consider the case where $(u_n)_n$ satisfies $ii)$. 
	
	\smallskip
	Set $v_n:=  F(u_n)$ for $n \in \na$ and $v:= F(u)$. 
	
	\smallskip
	By $ii)$,  $\int_{B_1}  e^{(4 \pi + t) u_n^2} \ dx$ is bounded for some $t>0$ and thus  $v_n$ is bounded in  $L^{s_0}(\R^2)$ with  
$s_0=1 + \frac{t}{4 \pi}>1$.
	Moreover, since 
	$$
	v_n \to v\qquad \text{in $L^1(\R^2)$,}
	$$
	interpolation yields that 
	$v_n \to v$ in $L^s(\re^2)$ for $1 \le s < s_0$.
	Finally	we conclude that 
	\begin{equation*}
		m_1(F) = \lim_{n \to \infty} \Phi(u_n) = \Phi(u)
	\end{equation*}
	so $m_1(F)$ is attained at $u \in \cB^*$.
		\end{proof}
	
	\bigskip
	
	Furthermore it is possible to derive qualitative property of the maximizer.

	\begin{lemma}
		Let $F$ satisfy $(A_0)$, and let $u \in \cB$ be a maximizer for $\Phi \big|_{\cB}$. Then, up to a change of sign, we have $u \in \cB^*$, and 
		$u$ is strictly positive in $B_1$.  
	\end{lemma}

	Moreover we can prove that the extremals are solutions of an Euler-Lagrange equation under some additional regularity properties on the function $F$.

	\begin{theorem}
		\label{e-l-theorem-section-ball-case}
		Suppose that $F \in C^1(\re)$ satisfies $(A_0)$, and suppose that $f:= F'$ satisfies
		\begin{equation*}
			f(t) \le c e^{\alpha t^2} \qquad \text{for $t \in \re$ with constants $c,\alpha>0$.}
		\end{equation*}
		Suppose furthermore that $u \in \cB$ is  a maximizer of $\Phi \big|_{\cB}$. Then there exists
		$\theta >0$
		such that $u$ satisfies the  Euler-Lagrange equation in weak sense, i.e., 
		$$
			\int_{B_1}\nabla u \nabla \phi \,dx  =  \theta \int_{B_1} \ln \frac{1}{|\cdot|}*\bigl(1_{B_1} F(u)\bigr) f(u)\phi\,dx
					$$
		for all $\phi \in H^1_0(B_1)$.
		Moreover, $u \in W^{2,p}(B_1) \cap C^{1,\sigma}(\overline {B_1})$ for all $p<\infty$, $\sigma \in (0,1)$. 
			\end{theorem}

	In the recent paper \cite{CiWeYu}, 
	we focus on the case	$F$ has at most $\beta=-1$ critical growth and
we turn to the question whether  the supremum  is attained in the case of critical growth. 
In some sense	we tackle the open problem if  there is a Log-Counterpart of Carleson-Chang paper.
\smallskip

We  recall that	a sequence $(u_n) \subset H^{1}_0(\Omega)$ is a normalized concentrating sequence if
	\begin{itemize}
		\item[i)] $\|\nabla u_n\|_2 =1$
		\item[ii)] $u_n \weak 0 \in H^1_0(\Omega)$
		\item[iii)] there exists $x_0 \in \Omega$ s.t. $\forall \rho >0, \int_{\Omega \setminus{B_\rho(x_0)}} |\nabla u_n|^2 dx \to 0$.
	\end{itemize}
We also remark that 	
$(u_n) \subset H^{1}_0(B_1)$ is a normalized concentrating sequence, we can suppose that $x_0=0$  and
	by symmetrization  $u_n=u_n^*$ is radially symmetric decreasing function.
	
	In \cite{CiWeYu} we give some answers to the above open problem, studying the behaviour of the nonlocal interaction functional on the concentrating sequences in the spirit of 
	the 
	celebrated paper
	\cite{figuereido-do-o-ruf}.
	
		\smallskip
		
		\begin{remark}
			In the case of unbounded domains,  the Trudinger-Moser type inequalities in \cite{CiWe2} are also proved for the entire domain $\re^2$ under an additional assumption
		\noindent
			\begin{itemize}
				\item[$(A_1)$] $F: \re \to [0,\infty)$ satisfies $(A_0)$, and $F(t)=O(|t|)$ as $t \to 0$.
			\end{itemize}	
		\end{remark}	
		
		\begin{remark}
		As written in \cite{CiWe2},	it is still an open  problem to derive Trudinger-Moser inequalities with Log-potentials for different geometry of bounded domains in critical cases. The  existence of maximizers may also  depend on the domain.
			It is still open the uniqueness, up to sign and translation, of maximizers of $\Phi|_{\cB}$.  
		\end{remark}
		
		\bigskip
		
\noindent
{\bf Acknowledgments.}

The authors are supported by PRIN PNRR  P2022YFAJH {\sl \lq\lq Linear and Nonlinear PDEs: New directions and applications''}, and partially supported by INdAM-GNAMPA.

The third author is supported by the INdAM-GNAMPA Project \lq \lq Fenomeni di concentrazione in PDEs non locali", CUP E53C22001930001.

The authors thank PNRR MUR project CN00000013 HUB - National Centre for HPC, Big Data and Quantum Computing (CUP H93C22000450007). 
\\


\begin{thebibliography}{6}

	\bibitem{AT} S. Adachi, K. Tanaka, Trudinger type inequalities in $\re^N$ and their best exponents,
				Proc. Amer. Math. Soc. 128,	(2000), {no. 7,}  2051-2057.
		
		
		
\bibitem{Adams} D. R. Adams,  A sharp inequality of J. Moser for higher order derivatives, Ann. of Math. 128, (1988), {no. 2,} 385–398.	
		
\bibitem{AD} A. Adimurthi, O. Druet, 
Blow-up analysis in dimension 2 and a sharp form of Moser–Trudinger inequality, Comm. Partial Differential Equations 29, (2004), {no. 1-2,} 295-322.
			
			
			
			
			
\bibitem{ASandeep}	A. Adimurthi, K. Sandeep,  A singular Moser–Trudinger embedding and its applications,  NoDEA Nonlinear Differential Equations Appl. 13, (2007), { no. 5-6,} 585-603.
			
			
			
		\bibitem{ABG} M. Agueh, S. Boroushaki, N. Ghoussoub, A dual Moser–Onofri inequality and its extensions to higher dimensional spheres, Ann. Fac. Sci. Toulouse Math. (6) 26, (2017), no. 2, 217–233.
		
		
\bibitem{Aubin} T. Aubin, Meilleures constantes dans le théorème  d'inclusion de Sobolev et un théorème  de Fredholm non
linéaire  pour la transformation conforme de la courbure scalaire, J. Funct. Anal. 32, (1979), { no. 2,} 148–174.
		

\bibitem{Bart}D. Bartolucci, Non-degeneracy, Mean Field Equations and the Onsager Theory of 2D Turbulence, Arch. Rational Mech. Anal. 230, (2018), 397–426.
		
		
		
\bibitem{BarMal} D. Bartolucci, A. Malchiodi, Mean field equations and domains of first kind, Rev. Mat. Iberoam. 38, (2022), {no. 4,}1067-1086.		
		
\bibitem{BatMan} L. Battaglia, G. Mancini, Remarks on the Moser–Trudinger inequality, Adv. Nonlinear Anal. 2, (2013), {no. 4,} 389–425. 

				
\bibitem{BVS} L. Battaglia, J. Van Schaftingen,	Groundstates of the Choquard equations with a sign-changing self-interaction potential, Z. Angew. Math. Phys. 69, (2018), {no. 3, 1-16.} 
		
		
		
	
			
\bibitem{Bek} W. Beckner, Sharp Sobolev inequalities on the sphere and the Moser-Trudinger inequality, Ann. of Math. (2) 138, (1993), no. 1, 213–242.
		
		
		
		
\bibitem{Berger} M. S. Berger, Riemannian structures of prescribed Gaussian curvature for compact 2-manifolds, J. Differential Geometry 5, (1971), 325-332.

		
		
\bibitem{BoCiSe} D. Bonheure, S. Cingolani, S. Secchi, Concentration phenomena for the Schr\"odinger-Poisson system in $\re^2$, Discrete Contin. Dyn. Syst. Ser. S 14, n. 5,(2021), 1631--1648.
	
			
\bibitem{BoCiVa} D. Bonheure, S. Cingolani, J. Van Schaftingen, The logarithmic Choquard equation: sharp asymptotics and nondegeneracy of the groundstate, J. Funct. Anal. 272, (2017), 5255-5281.
		
		
	
	
\bibitem{BCM}	N. Borgia, S. Cingolani, G. Mancini, On $N$-Euclidean Logarithmic Moser-Trudinger-Onofri inequality and some geometrical variants,   in preparation.
		

\bibitem{CLMP0} E. Caglioti, P. L. Lions, C. Marchioro, M.	Pulvirenti, A special class of stationary flows for two-dimensional euler equations: A statistical mechanics description., Comm. Math. Phys. 143, (1992), 501-525.
		
\bibitem{CLMP} E. Caglioti, P. L. Lions, C. Marchioro, M.	Pulvirenti, A special class of stationary flows for two-dimensional euler equations: A statistical mechanics description. Part II, Comm. Math. Phys. 174, (1995), 229-260.


\bibitem{CaRuf} M. Calanchi, B. Ruf, On Trudinger–Moser type inequalities with logarithmic weights, J. Differential Equations 258, n. 6, (2015), 1967-1989.
	
\bibitem{cao} {D. M. Cao}, Nontrivial solution of semilinear elliptic equation with critical exponent in $\re^2$, Comm. P.D.E. 17, (1992), 407-435.
	
\bibitem{carlenloss} E. Carlen, M. Loss, Competing symmetries, the logarithmic HLS inequality and Onofri's inequality on $S^n$, Geometric \& Functional Analysis GAFA, 1992 - Springer. 
	
	
	
	
\bibitem{carlesonchang} L. Carleson, A. Chang,	On the existence of an extremal function for an inequality of J. Moser, Bull. Sci. Math. 110, (1986), 113-127.
	
	

	
\bibitem{CST} D. Cassani, F. Sani, C. Tarsi, Equivalent Moser type inequalities in $\re^2$ and the zero mass case, J. Funct. Anal.  267, (2014), 4236-4263.
	
	
	
	
\bibitem{CY1} S.-Y. A. Chang , P. C. Yang, 
Prescribing Gaussian curvature on $S^2$, Acta Mathematica 159 (1), (1987), 215-259.

\bibitem{CY2} S.-Y. A. Chang , P. C. Yang, Conformal deformation of metrics on $S^2$, Journal of Differential Geometry 27 (2), (1988), 259-296.

\bibitem{CCL} S.-Y. A. Chang, C.-C. Chen, C.-S. Lin, Extremal functions for a mean field equation in two dimension, Lectures on Partial Differential Equations in Honor of Louis Nirenberg's 75th Birthday, Internat. Press (2003)

\bibitem{ChenSing} W. Chen, A trudinger inequality on surfaces with conical singularities, Proceedings of the American Mathematical Society, 108 (3), (1990), 821-832.

\bibitem{CL} W. Chen , C. Li,  Prescribing gaussian curvatures on surfaces with conical singularities. J Geom Anal 1, (1991), 359–372.

\bibitem{CL1} C.-C. Chen, C.-S. Lin, Sharp Estimates for Solutions of Multi-Bubbles in Compact Riemann Surfaces, Communications on Pure and Applied Mathematics, 55 (6), (2002), 728-771.


\bibitem{CL2} C.-C. Chen, C.-S. Lin, Topological Degree for a Mean Field Equation on Riemann Surfaces, Communications on Pure and Applied Mathematics, 56 (12), (2003) , 1667-1727.

	
	
\bibitem{CiJe} S. Cingolani, L. Jeanjean, Stationary solutions with prescribed $L^2$-norm for the planar Schr\"odinger - Poisson system, SIAM J. Math. Anal. 51, (2019), 3533-3568.
	
	
\bibitem{CiWe} S. Cingolani,  T. Weth, On the planar Schr\"odinger-Poisson system, Ann. Inst. H. Poincar\'e Anal. Non Lin\'eaire, 33, (2016), 169-197.
	
	
\bibitem{CiWe2}S. Cingolani, T. Weth, Trudinger-Moser-type inequality with logarithmic convolution potentials, J. London Mathematical Society, 105, (2022), no. 3, 1897-1935.
	
	

\bibitem{CiWeYu} S. Cingolani, T. Weth, M. Yu, Extremal functions for the critical Trudinger-Moser inequality with logarithmic kernels, submitted paper.
	
	
	
\bibitem{CR}	G. Csat\'o, P. Roy, Extremal functions for the singular Moser–Trudinger inequality in 2 dimensions, Calc. Var. Partial Differential Equations 54, (2015),  2341-2366.
	
	
\bibitem{figuereido-do-o-ruf} D.G. de Figueiredo, M. do O, B. Ruf,  On an inequality by N. Trudinger and J. Moser and related elliptic equations, Communications on Pure and Applied Mathematics 55, (2002), 135-152.	
	
\bibitem{DD} M. del Pino, J. Dolbeault, The Euclidean Onofri inequality in higher dimensions Int. Math. Res. Not. IMRN 2013 (2013), no. 15, 3600–3611.
			
\bibitem{ManDLT} A. DelaTorre, G. Mancini, Improved Adams-type inequalities and their extremals in dimension 2m, Commun. Contemp. Math. 23, no. 05, (2021), 2050043. 			
			
\bibitem{doO} J.M. do \'O, N-Laplacian equations in with critical growth, Abstr. Appl. Anal. 2, (1997), 301-315.
			
			
\bibitem{DEJ} J. Dolbeault, M.J. Esteban, G. Jankowiak, The Moser-Trudinger-Onofri Inequality, Chinese Ann. Math. Ser. B36, (2015), no. 5,  777–802.
			
	
	
\bibitem{Dolbeault-Perthame} J. Dolbeault, B. Perthame, Optimal critical mass in the two dimensional Keller-Segel model in $\re^2$. C. R. Acad. Sci. Paris, Ser. I 339, (2004), 611-616.
	


\bibitem{DongLamLu} M. Dong, N. Lam, G. Lu, Sharp weighted Trudinger–Moser and Caffarelli–Kohn–Nirenberg inequalities and their extremal functions, Nonlinear Analysis 173, (2018), 75-98.
	
	
	
\bibitem{DuWe} M. Du, T. Weth, Ground states and high energy solutions of the planar Schr\"odinger-Poisson system, Nonlinearity 30, (2017), 3492-3515.


	

	
		
\bibitem{Fontana}L. Fontana, Sharp borderline Sobolev inequalities on compact  Riemannian manifolds, Comment. Math. Helv. 68, (1993), 415-454.
	
	
\bibitem{flucher} M. Flucher, Extremal functions for the Trudinger-Moser inequality in 2 dimensions, Comm. Math. Helv.	67, (1992), 471-479.
	
	
\bibitem{IMN} S. Ibrahim, N. Masmoudi, K. Nakanishi, Trudinger–Moser inequality on the whole plane with the exact growth condition, J. Eur. Math. Soc. 17, (2015), 819–835. 

\bibitem{IMNS}S. Ibrahim, N. Masmoudi, K. Nakanishi, F. Sani, Sharp threshold nonlinearity for maximizing the Trudinger-Moser inequalities, J. Funct. Anal.278 (2020), no. 1, 108302, 52 pp. 
	
	
\bibitem{Ishi} M. Ishiwata, Existence and nonexistence of maximizers for variational problems associated with Trudinger–Moser type inequalities in $\R^N$, Math. Ann. 351 (2011), no. 4, 781–804.	
	
\bibitem{IM} S. Iula,  G. Mancini, Extremal functions for singular Moser–Trudinger embeddings, Nonlinear Analysis 156, (2017), 215-248.
		
\bibitem{JM} G. Joyce, D. Montgomery, Negative temperature states for the two-dimensional guiding-cente plasma, J. Plasma Phys. 10, (1973), 107-121.		
		
		
\bibitem{KW}J. L. Kazdan, F. W. Warner, Curvature Functions for Compact 2-Manifolds, Annals of Mathematics, Second Series, Vol. 99,  (1974), no. 1, 14-47.
		
		
\bibitem{LamLu}	N. Lam, G. Lu, Weighted Moser–Onofri–Beckner and Logarithmic Sobolev Inequalities, J. Geometrical Analysis, J Geom Anal.  28, (2018), 1687–1715.
		
		
\bibitem{Li1} Y. X.	Li,	Moser-Trudinger inequality on compact Riemannian manifolds of dimension two, J. Partial Differential Equations  14,  (2001), no. 2, 163-192.

\bibitem{Li2} Y.X. Li, Extremal functions for the Moser-Trudinger inequalities on compact Riemannian manifolds, Sci. China Ser. A 48 (5), (2005), 618–648.


\bibitem{LiLu} J.Li, G.Lu, Critical and subcritical Trudinger-Moser inequalities on complete noncompact Riemannian manifolds,
Advances in Mathematics  389, (2021), no. 8, 107915.
		
		
		
		
		
\bibitem{LiRuf}	Y.Li, B.Ruf, A sharp Trudinger-Moser type inequality for unbounded domains in $\re^N$, Indiana University Mathematics Journal,{57, (2008), no. 1, 451–480.}

\bibitem{LiYang}	X. Li, Y. Yang Extremal functions for singular Trudinger–Moser inequalities in the entire Euclidean space J. Differential Equations 264, (2018), no. 8, 4901-4943.		
		
\bibitem{Lions} P.L. Lions, The concentration-compactness principle in the calculus of variations. The limit case, part 1, Riv. Mat. Iberoamericana 1, (1985), 145-201.
			
			
\bibitem{YoungAdams}	G. Lu, Y. Yang, Adams' inequalities for bi-Laplacian and extremal functions in dimension four, Adv. Math. 220, (2009), 1135-1170.
	
\bibitem{Malchiodi} A. Malchiodi, Topological methods for an elliptic equation with exponential nonlinearities,
Discrete Contin. Dyn. Syst. 21, no. 1,  (2008), 277-294.

\bibitem{MalRuiz} A. Malchiodi, D. Ruiz, New Improved Moser–Trudinger Inequalities and Singular Liouville Equations on Compact Surfaces,  Geom. Funct. Anal. Vol. 21, (2011), 1196–1217.
	
	
\bibitem{MM} G. Mancini, L. Martinazzi, The Moser-Trudinger inequality and its extremals on a disk via energy estimates, Calc. Var. {56, (2017), no. 4, Paper No. 94, 26 pp.} 



\bibitem{MM2} G. Mancini, L. Martinazzi, Extremals for Fractional Moser–Trudinger Inequalities in Dimension 1 via Harmonic Extensions and Commutator Estimates, Advanced Nonlinear Studies 20, (2020), no. 3, 599-632.

\bibitem{ManSand} G. Mancini,  K. Sandeep, Moser–Trudinger inequality on conformal discs, Commun.
Contemp. Math. 12, (2010), no. 6, 1055–1068.	
	
\bibitem{ManThi} G. Mancini, P.-D. Thizy, Non-existence of extremals for the Adimurthi–Druet inequality, Journal of Differential Equations {266, (2019), no. 2–3, 1051-1072.}
	
\bibitem{Martinazzi} L. Martinazzi, Fractional Adams–Moser–Trudinger type inequalities, Nonlinear Anal. 127, (2015), 263–278.		
		
\bibitem{Masaki} S. Masaki, Local existence and WKB approximation of solutions to Schr\"odinger-Poisson system in the two-dimensional whole space, Communications P.D.E. 35, (2010), no. 12, 2253-2278.
		
		
\bibitem{Masaki2} S. Masaki, Energy Solution to a Schr\"odinger-Poisson System in the Two-Dimensional Whole Space, SIAM J. Math. Anal. 43, (2011), no. 6, 2719-2731.
		
		
	
		
			
\bibitem{M} J. Moser, A sharp form of an inequality by N. Trudinger, Indiana Univ. Math. J. 20, (1971), no. 11, 1077-1092.
			
\bibitem{M2} J. Moser, On a nonlinear problem in differential geometry, Dynamical systems. Academic Press, (1973), 273-280.
			
\bibitem{O} E. Onofri, On the positivity of the effective action in a theory of random surfaces, Comm. Math. Phys. 86, (1982), no. 3, 321–326.

\bibitem{OPS}
B.~Osgood, R.~Phillips,  P.~Sarnak,  Extremals of determinants of {L}aplacians,  J. Funct. Anal. 80, (1988), no. 1, 148-211.

\bibitem{OPS2}
B.~Osgood, R.~Phillips,  P.~Sarnak, Compact isospectral sets of surfaces, J. Funct. Anal. 80, (1988), no. 1, 212-234.			
			
\bibitem{poho}{S. I. Pohozaev, On the eigenfunctions of the equation $\Delta u + \lambda f(u)=0$, 
		Dokl. Akad. Nauk SSSR 165, (1965), 36-39 (in Russian).	}		
			
\bibitem{Poly1}
A. M. Polyakov, Quantum geometry of bosonic strings, Phys. Lett. B 103, (1981), no. 3, 207-210.

\bibitem{Poly2}
A.~M. Polyakov, Quantum geometry of fermionic strings, Phys. Lett. B 103, (1981), no. 3, 211-213.
			
	
	
	
\bibitem{pruss} A. R. Pruss, Nonexistence of maxima for perturbations of some inequalities with critical growth, Can.
Math. Bull. 39, (1996), 227–237. 		
	
		\bibitem{ruf} B. Ruf, A sharp Trudinger-Moser type inequality for unbounded domains in $\re^2$,  J. Funct. Anal. 219, (2005), {no. 2,} 340-367.
	
\bibitem{Struwe} M. Struwe, Critical points of embeddings of $H^{1,n}_0$ into Orlicz spaces,  Ann. Inst. H. Poincar\'e Anal. Non Lineaire 5, (1988), {no. 5,} 425-464.
	
\bibitem{StruweCurv} M. Struwe, A flow approach to Nirenberg’s problem, Duke Math. J. 128, (2005), no. 1,  19-64.
	
	
\bibitem{stubbe} J. Stubbe, Bound states of two-dimensional Schr\"odinger-Newton equations, ArXive 0807.4059v1 (2008)
	
	
	
\bibitem{SV} J. Stubbe, M. Vuffray, Bound states of the  Schr\"odinger- Newton model in low dimensions, Nonlinear Analysis 73, (2010), {no. 10,} 3171-3178.

	
\bibitem{Suzuki} T. Suzuki,
Free energy and self-interacting particles, Progress in Nonlinear Differential Equations and their Applications, Birkh\"auser, Boston 62, (2005).
	
	
\bibitem{Thizy} P. D. Thizy, When does a perturbed Moser-Trudinger inequality admit an extremal?, Anal. PDE . 13, (2020), no. 5, 1371–1415.


\bibitem{Tintarev} C. Tintarev,
Trudinger–Moser inequality with remainder terms, J. Funct. Anal. 266, (2014), {no. 1,} 55-66.
	
\bibitem{Troy} M. Troyanov, Prescribing curvature on compact surfaces with conical singularities, Trans. Amer. Math. Soc. 324, (1991), {no. 2,} 793-821.
	
	
\bibitem{Tru} N.S. Trudinger, On imbeddings into Orlicz spaces and some applications, J. Math. Mech. 17, (1967), no. 5, 473-483.
	
	
	
\bibitem{W} G. Wolansky, On steady distributions of self-attracting clusters under friction and fluctuations, Arch. Rational Mech. Anal. 119, (1992), no. 4, 355-391.
	

	
	
\bibitem{Yang}	Y. Yang, Extremal functions for Trudinger–Moser inequalities of Adimurthi–Druet type in dimension two, J. Differential Equations 258, (2015), {no. 9,} 3161-3193.
	
	
	

	
\bibitem{YangZhu}	Y. Yang, X. Zhu, Blow-up analysis concerning singular Trudinger–Moser inequalities in dimension two,
	J. Funct. Anal. 272, (2017), {no. 8,}  3347-3374.
	
	
\bibitem{Yudo} V. I. Yudovich,  Some estimates connected with integral operators and with solutions of elliptic equations,
Dokl. Akad. Nauk SSSR 138,  (1961),{no. 4,} 805–808 (in Russian).
	
	
	
	
	
			
			
\end{thebibliography}
\end{document}